\documentclass[final,1p,times]{elsarticle}

\usepackage{amssymb}
\usepackage{amsmath}

\journal{Nonlinear Analysis, Theory, Methods \& Applications}

\newtheorem{theorem}{Theorem}
\newtheorem{proposition}[theorem]{Proposition}
\newtheorem{lemma}[theorem]{Lemma}
\newtheorem{corollary}[theorem]{Corollary}

\newdefinition{definition}{Definition}
\newdefinition{remark}{Remark}

\newproof{proof}{Proof}

\begin{document}

\begin{frontmatter}

\title{Calculus of Variations on Time Scales
with Nabla Derivatives\tnoteref{ack}}

\tnotetext[ack]{Supported by the Portuguese Foundation
for Science and Technology (FCT), through the R\&D unit
\emph{Centre for Research on Optimization and Control} (CEOC),
cofinanced by the European Community Fund FEDER/POCI 2010.}

\author{Nat\'{a}lia Martins}
\ead{natalia@ua.pt}

\author{Delfim F. M. Torres}
\ead{delfim@ua.pt}

\address{Department of Mathematics,
University of Aveiro,
3810-193 Aveiro, Portugal}


\begin{abstract}
We prove a necessary optimality condition of Euler-Lagrange type
for variational problems on time scales involving nabla
derivatives of higher-order. The proof is done using a new and
more general fundamental lemma of the calculus of variations on
time scales.
\end{abstract}


\end{frontmatter}


\section{Introduction}

The theory of \emph{time scales} was born in 1988 with the
work of Stephan Hilger \cite{Hilger2}, providing a rich
theory that unify and extend discrete and continuous analysis
\cite{Bohner-Peterson1,Bohner-Peterson2}.
The study of the calculus of variations in the context of time scales
has it's beginning in 2004 with the paper \cite{Bohner:2004} of
Martin Bohner. The main result of \cite{Bohner:2004}
is an Euler-Lagrange necessary optimality equation for a first order
variational problem involving delta derivatives on a time scale $\mathbb{T}$:
\begin{theorem}[\cite{Bohner:2004}]
\label{Euler-Lagrange: Bohner} If $y_\ast \in C_{\textrm{rd}}^2$
is a weak local extremum  of the problem
\begin{equation}
\label{eq:prb:cv:delta}
\mathcal{L}[y(\cdot)]
=\int_{a}^{b}L(t,y^\sigma(t),y^\Delta(t))\Delta t
\longrightarrow  \mathrm{extr} \, , \quad
y(a)=\alpha, \ y(b)=\beta,
\end{equation}
where $(t,u,v) \rightarrow L(t,u,v)$ is a $C^2$ function,
then the Euler-Lagrange equation
\begin{equation*}
L_{y^\Delta}^\Delta(t,y^\sigma_\ast(t),y_\ast^\Delta(t))
=L_{y^\sigma}(t,y^\sigma_\ast(t),y_\ast^\Delta(t))
\end{equation*}
holds for $t\in[a,b]^{\kappa}$.
\end{theorem}
Since the pioneer work \cite{Bohner:2004} of Martin Bohner,
Theorem~\ref{Euler-Lagrange: Bohner} has been extended
in several different directions in order
to analyze variational problems on time scales with:
(i) non-fixed boundary conditions \cite{zeidan:2004};
(ii) two independent variables \cite{BonGus};
(iii) higher-order delta derivatives \cite{Lisboa07};
(iv) an invariant group of parameter-transformations \cite{Zbig};
(v) multiobjectives \cite{Aveiro2Basia};
(vi) isoperimetric constraints \cite{iso:ts}.
A different direction of study, which seems of special interest
to applications, in particular to economics, is given in \cite{Atici:2006}
where a Euler-Lagrange type equation is obtained
for a first order variational problem on time scales
involving nabla derivatives instead of delta ones:
\begin{theorem}[\cite{Atici:2006}]
\label{thm:cv:nabla}
If a function $y_\ast \in C^2$ provides a weak local extremum
of the problem
\begin{equation}
\label{eq:prb:cv:nabla}
\mathcal{L}[y(\cdot)]=\int_{\rho^2(a)}^{\rho(b)}L(t,y^\rho(t),y^\nabla(t))\nabla
t \longrightarrow  \mathrm{extr} \, , \quad
y(\rho^2(a))=\alpha \, , \  y(\rho(b))=\beta,
\end{equation}
where $(t,u,v) \rightarrow L(t,u,v)$ is a $C^2$ function
of $(u,v)$ for each $t \in [\rho^2(a),\rho(b)] \subseteq
\mathbb{T}$, then $y_\ast$ satisfies the
Euler-Lagrange equation
$$
L_{y^\nabla}^\nabla(t,y^\rho_\ast(t),y_\ast^\nabla(t))
=L_{y^\rho}(t,y^\rho_\ast(t),y_\ast^\nabla(t))
$$
for $t \in [\rho(a), b]$.
\end{theorem}
The proof of Theorem~\ref{thm:cv:nabla} found in \cite{Atici:2006}
has, however, some inconsistences \cite{Tangier07}. Moreover, the
nabla problem \eqref{eq:prb:cv:nabla} is defined in a way that is
not completely analogous to the more studied and established delta
problem \eqref{eq:prb:cv:delta} (\textrm{i.e.}, one does not
obtain \eqref{eq:prb:cv:nabla} from \eqref{eq:prb:cv:delta} by
simply substituting the forward jump operator $\sigma$ by the
backward jump operator $\rho$ and the delta derivative $\Delta$ by
the nabla derivative $\nabla$). The main goal of this paper is to
generalize the results of \cite{Atici:2006} to higher-order nabla
variational problems on time scales in a consistent way with the
delta theory. This is done by first proving some new fundamental
lemmas of the calculus of variations on time scales that fix the
inconsistencies of \cite{Atici:2006} pointed out in
\cite{Tangier07}. Compared with the delta approach followed in
\cite{Lisboa07}, our technique provides a simpler and more direct
proof to the higher-order Euler-Lagrange equations. Moreover, our
proof seems to be new even for the continuous time $\mathbb{T} =
\mathbb{R}$.

The paper is organized as follows. In Section~\ref{sec:pr} we collect
all the necessary elements of the nabla calculus on time scales.
In Section~\ref{sec:mr} we state and prove our results:
in \S\ref{sub:sec:fl} we prove a new and more
general fundamental lemma of the calculus of variations on time
scales (Lemma~\ref{lemma7}); in \S\ref{sub:sec:EL:ho} we obtain a
higher-order nabla differential Euler-Lagrange equation
(Theorem~\ref{Euler-Lagrange for higher order}).
As an example, we give
the Euler-Lagrange equation for the q-calculus variational problem
(Corollary~\ref{sub:sec:example}).


\section{Preliminary results}
\label{sec:pr}

For a general introduction to the calculus on time scales we refer
the reader to the books \cite{Bohner-Peterson1,Bohner-Peterson2}.
Here we only give those notions and results needed in the sequel.
More precisely, we are interested in the nabla approach to time scales
\cite{AticiGuseinov:2002}. As usual,  $\mathbb{R}$, $\mathbb{Z}$,
and $\mathbb{N}$ denote, respectively, the set of real,
integer, and natural numbers.

A \emph{Time Scale} $\mathbb{T}$ is an arbitrary non empty closed
subset of $\mathbb{R}$. Thus, $\mathbb{R}$, $\mathbb{Z}$, and
$\mathbb{N}$, are trivial examples of times scales. Other examples
of times scales are: $[-1,4] \bigcup \mathbb{N}$,
$h\mathbb{Z}:=\{h z | z \in \mathbb{Z}\}$ for some $h>0$,
$q^{\mathbb{N}_0}:=\{q^k | k \in \mathbb{N}_0\}$ for some $q>1$,
and the Cantor set. We assume that a time scale $\mathbb{T}$ has
the topology that it inherits from the real numbers with the
standard topology.

The \emph{forward jump operator}
$\sigma:\mathbb{T}\rightarrow\mathbb{T}$ is defined by
$\sigma(t)=\inf{\{s\in\mathbb{T}:s>t}\}$
if $t\neq \sup \mathbb{T}$, and $\sigma(\sup \mathbb{T})=\sup
\mathbb{T}.$ The \emph{backward jump operator}
$\rho:\mathbb{T}\rightarrow\mathbb{T}$ is defined by
$\rho(t)=\sup{\{s\in\mathbb{T}:s<t}\}$
if $t\neq \inf \mathbb{T}$, and $\rho(\inf \mathbb{T})=\inf
\mathbb{T}.$

A point $t\in\mathbb{T}$ is called \emph{right-dense},
\emph{right-scattered}, \emph{left-dense} and
\emph{left-scattered} if $\sigma(t)=t$, $\sigma(t)>t$, $\rho(t)=t$,
and $\rho(t)<t$, respectively. We say that $t$ is \emph{isolated}
if $\rho(t)<t<\sigma(t)$, that $t$ is \emph{dense} if $\rho(t)=t=\sigma(t)$.
The \emph{(backward) graininess function}
$\nu:\mathbb{T}\rightarrow[0,\infty)$ is defined by
$\nu(t)=t - \rho(t)$, for all $t\in\mathbb{T}$.
Hence, for a given $t$, $\nu(t)$ measures the distance
of $t$ to its left neighbor.
It is clear that when $\mathbb{T}=\mathbb{R}$ one has
$\sigma(t)=t=\rho(t)$, and $\nu(t)=0$
for any $t$. When $\mathbb{T}=\mathbb{Z}$,
$\sigma(t)=t+1$, $\rho(t)=t-1$, and $\nu(t)=1$
for any $t$.

In order to introduce the definition of nabla derivative, we
define a new set $\mathbb{T}_\kappa$ which is derived from $\mathbb{T}$
as follows: if  $\mathbb{T}$ has a right-scattered minimum $m$,
then $\mathbb{T}_\kappa=\mathbb{T}\setminus\{m\}$; otherwise,
$\mathbb{T}_\kappa= \mathbb{T}.$

\begin{definition}
We say that a function $f:\mathbb{T}\rightarrow\mathbb{R}$ is
\emph{nabla differentiable} at $t\in\mathbb{T}_\kappa$ if there is a
number $f^{\nabla}(t)$ such that for all $\varepsilon>0$ there
exists a neighborhood $U$ of $t$ (\textrm{i.e.},
$U=]t-\delta,t+\delta[\cap\mathbb{T}$ for some $\delta>0$) such
that
$$|f(\rho(t))-f(s)-f^{\nabla}(t)(\rho(t)-s)|
\leq\varepsilon|\rho(t)-s|,\mbox{ for all $s\in U$}.$$ We call
$f^{\nabla}(t)$ the \emph{nabla derivative} of $f$ at $t$.
Moreover, we say that $f$ is \emph{nabla differentiable} on
$\mathbb{T}$ provided $f^{\nabla}(t)$ exists for all $t \in
\mathbb{T}_\kappa$.
\end{definition}

\begin{theorem}
\label{propriedades derivada}
Let $\mathbb{T}$ be a time scale, $f:\mathbb{T}\rightarrow\mathbb{R}$,
and $t\in\mathbb{T}_\kappa$. The following holds:
\begin{enumerate}
\item If $f$ is nabla differentiable at $t$, then $f$ is
    continuous at $t$.
\item If $f$ is continuous at $t$ and $t$ is left-scattered,
    then $f$ is nabla differentiable at $t$ and
$$f^{\nabla}(t)=\frac{f(t)-f(\rho(t))}{t-\rho(t)}.$$

\item If $t$ is left-dense, then $f$ is nabla differentiable
    at $t$ if and only if the limit
$$
\lim_{s\rightarrow t} \frac{f(t)-f(s)}{t-s}
$$
exists as a finite number. In this case,
$$
f^\nabla(t)=\lim_{s\rightarrow t} \frac{f(t)-f(s)}{t-s}.
$$
\item If $f$ is nabla differentiable at $t$, then
$f(\rho(t))=f(t)-\nu(t)f^\nabla(t)$.
\end{enumerate}
\end{theorem}

\begin{remark}
When $\mathbb{T}=\mathbb{R}$, then $f:\mathbb{R}
\rightarrow \mathbb{R}$ is nabla differentiable  at $t \in
\mathbb{R}$ if and only if $$ \displaystyle f^{\nabla}(t)=
\lim_{s\rightarrow t}\frac{f(t)-f(s)}{t-s}$$ exists, \textrm{i.e.},
if and only if $f$ is differentiable at $t$ in the ordinary
sense. When $\mathbb{T}=\mathbb{Z}$, then $f:\mathbb{Z}
\rightarrow \mathbb{R}$ is always nabla differentiable at
$t \in \mathbb{Z}$ and
$$
f^{\nabla}(t)=\frac{f(t)-f(\rho(t))}{t-\rho(t)}=
f(t)-f(t-1),
$$
\textrm{i.e.}, $f^{\nabla}$ is the usual backward difference.
For any time scale $\mathbb{T}$,
when $f$ is a constant, then $f^{\nabla}=0$;
if $f(t)=k t$ for some constant $k$, then $f^{\nabla}=k$.
\end{remark}

In order to simplify expressions, we denote the composition
$f\circ \rho$ by $f^{\rho}$. We also use the standard
notation $\rho^0(t)=t$, $\sigma^0(t)=t$, $\rho^n(t)=(\rho \circ
\rho^{n-1})(t)$, and $\sigma^n(t)=(\sigma \circ \sigma^{n-1})(t)$
for $n \in \mathbb{N}$.

\begin{theorem}
Suppose $f,g:\mathbb{T}\rightarrow\mathbb{R}$ are
nabla differentiable at $t\in\mathbb{T}_\kappa$. Then,
\begin{enumerate}
\item the sum $f+g:\mathbb{T}\rightarrow\mathbb{R}$ is nabla
differentiable at $t$ and
$(f+g)^{\nabla}(t)=f^{\nabla}(t) + g^{\nabla}(t)$;

\item for any constant $\alpha$, $\alpha
f:\mathbb{T}\rightarrow\mathbb{R}$ is nabla differentiable
at $t$ and $(\alpha f)^{\nabla} (t)=\alpha f^{\nabla}(t)$;

\item the product $fg:\mathbb{T}\rightarrow\mathbb{R}$ is
    nabla differentiable at $t$ and
$$
\begin{array}{rcl}
(fg)^{\nabla}(t)& = & f^{\nabla}(t)g(t) +
f^{\rho}(t)g^{\nabla}(t)\\ &&\\ &=& f^{\nabla}(t)g^{\rho}(t)+
f(t)g^{\nabla}(t).
\end{array}
$$
\end{enumerate}
\end{theorem}

Nabla derivatives of higher order are defined in
the standard way: we define the $r^{th}-$\emph{nabla derivative}
($r\in\mathbb{N}$) of $f$ to be the function
$f^{\nabla^r}:\mathbb{T}_{\kappa^r}\rightarrow\mathbb{R}$, provided
$f^{\nabla^{r-1}}$ is nabla differentiable on
$\mathbb{T}_{\kappa^r}:=\left(\mathbb{T}_{\kappa^{r-1}}\right)_\kappa$.

\begin{definition}
A function $F:\mathbb{T}\rightarrow\mathbb{R}$ is called a
\emph{nabla antiderivative} of
$f:\mathbb{T}\rightarrow\mathbb{R}$ provided
$$
F^{\nabla}(t)=f(t), \qquad  \forall t \in \mathbb{T}_\kappa.
$$
In this case we define the \emph{nabla integral} of $f$ from $a$
to $b$ ($a,b \in \mathbb{T}$) by
$$\int_{a}^{b}f(t)\nabla t:=F(b)-F(a).$$
\end{definition}

In order to present a class of functions that possess a nabla
antiderivative, the following definition is introduced:

\begin{definition}
Let $\mathbb{T}$ be a time scale, $f:\mathbb{T}\rightarrow\mathbb{R}$.
We say that function $f$ is \emph{ld-continuous} if it is continuous
at the left-dense points and its right-sided limits exist (finite) at all
right-dense points.
\end{definition}

Some results concerning ld-continuity are useful:

\begin{theorem}
Let $\mathbb{T}$ be a time scale,
$f:\mathbb{T}\rightarrow\mathbb{R}$.
\begin{enumerate}
\item If $f$ is continuous, then $f$ is ld-continuous.

\item The backward jump operator $\rho$ is ld-continuous.

\item If $f$ is ld-continuous, then $f^{\rho}$ is also
ld-continuous.

\item If $\mathbb{T}=\mathbb{R}$, then $f$ is continuous if
and only if $f$ is ld-continuous.

\item If $\mathbb{T}=\mathbb{Z}$, then $f$ is ld-continuous.
\end{enumerate}
\end{theorem}

\begin{theorem}
\label{antiderivada}
Every ld-continuous function has a nabla
antiderivative. In particular, if $a \in \mathbb{T}$, then the
function $F$ defined by
$$
F(t)= \int_{a}^{t}f(\tau)\nabla\tau, \quad t \in \mathbb{T} \, ,
$$
is a nabla antiderivative of $f$.
\end{theorem}

The set of all ld-continuous functions
$f:\mathbb{T}\rightarrow\mathbb{R}$ is denoted by
$C_{\textrm{ld}}(\mathbb{T}, \mathbb{R})$, and the set of all nabla
differentiable functions with ld-continuous derivative by
$C_{\textrm{ld}}^1(\mathbb{T}, \mathbb{R})$.

\begin{theorem}
\label{propriedade integral}
If $f \in C_{\textrm{ld}}(\mathbb{T}, \mathbb{R})$
and $t \in \mathbb{T}_\kappa$, then
$$
\int_{\rho(t)}^{t}f(\tau)\nabla\tau=\nu(t)f(t).
$$
\end{theorem}

\begin{theorem}
\label{propriedades nabla integral}
If $a,b,c \in \mathbb{T}$, $a \le c \le b$, $\alpha \in \mathbb{R}$,
and $f,g \in C_{\textrm{ld}}(\mathbb{T}, \mathbb{R})$, then
\begin{enumerate}
\item $\displaystyle \int_{a}^{b}\left(f(t) + g(t) \right)
    \nabla t= \int_{a}^{b}f(t)\nabla t +
    \int_{a}^{b}g(t)\nabla t$;

\item $\displaystyle \int_{a}^{b} \alpha f(t)\nabla t =\alpha
    \int_{a}^{b}f(t)\nabla t$;

\item $\displaystyle \int_{a}^{b}  f(t)\nabla t = -
    \int_{b}^{a} f(t)\nabla t$;

\item $\displaystyle \int_{a}^{a}  f(t)\nabla t=0$;

\item $\displaystyle \int_{a}^{b}  f(t)\nabla t =
    \int_{a}^{c}  f(t)\nabla t + \int_{c}^{b} f(t)\nabla t$;

\item If $f(t)> 0$ for all $a < t \leq b$, then $
    \displaystyle \int_{a}^{b}  f(t)\nabla t > 0$;

\item $\displaystyle \int_{a}^{b}f^\rho(t)g^{\nabla}(t)\nabla t
=\left[(fg)(t)\right]_{t=a}^{t=b}-\int_{a}^{b}f^{\nabla}(t)g(t)\nabla t$;

\item $\displaystyle \int_{a}^{b}f(t)g^{\nabla}(t)\nabla t
=\left[(fg)(t)\right]_{t=a}^{t=b}-\int_{a}^{b}f^{\nabla}(t)g^\rho(t)\nabla t$.
\end{enumerate}
\end{theorem}

The last two formulas on Theorem~\ref{propriedades nabla integral}
are usually called \emph{integration by parts formulas}. These formulas
are used several times in this work.

\begin{remark}
Let $a, b \in \mathbb{T}$ and $f\in
C_{ld}(\mathbb{T}, \mathbb{R})$.
For $\mathbb{T}=\mathbb{R}$, then
$\int_{a}^{b}f(t)\nabla t = \int_{a}^{b}f(t) dt$, where
the integral on the right side is the usual Riemann
integral. For $\mathbb{T}=\mathbb{Z}$, then $\displaystyle
\int_{a}^{b}f(t)\nabla t = \sum_{t=a+1}^{b}f(t)$ if $a<b$,
$\displaystyle \int_{a}^{b}f(t)\nabla t=0$ if $a=b$, and
$\displaystyle \int_{a}^{b}f(t)\nabla t
= - \sum_{t=b+1}^{a}f(t)$ if $a>b$.
\end{remark}

Let $a, b \in \mathbb{T}$ with $a<b$. We define the interval
$[a,b]$ in $\mathbb{T}$ by
$$
[a,b]:=\{t \in \mathbb{T}: a \leq t
\leq b\}.
$$
Open intervals and half-open intervals in $\mathbb{T}$ are defined
accordingly. Note that $[a,b]_\kappa=[a,b]$ if $a$ is right-dense and
$[a,b]_k=[\sigma(a),b]$ if $a$ is right-scattered.


\section{Main results}
\label{sec:mr}

Our main objective is to establish a
necessary optimality condition for problems of the
calculus of variations with higher-order nabla derivatives.
We formulate the higher-order variational problem with nabla derivatives
as follows:
\begin{equation}
\label{problema varacional}
\begin{gathered}
\mathcal{L}[y(\cdot)]=\int_{\sigma^{r-1}(a)}^{b}
L(t,y^{\rho^r}(t),y^{\rho^{r-1}\nabla}(t),\ldots,
y^{\rho\nabla^{r-1}}(t),y^{\nabla^r}(t))\nabla t\longrightarrow
\mathrm{extr} \, ,\\
y(\sigma^{r-1}(a))=\alpha_0,  \quad  y(b)=\beta_0, \\
\vdots\\ \tag{P} y^{\nabla^{r-1}}(\sigma^{r-1}(a))=\alpha_{r-1},
\quad y^{\nabla^{r-1}}(b)=\beta_{r-1},
\end{gathered}
\end{equation}
where $r\in\mathbb{N}$. We assume that:
\begin{enumerate}

\item The admissible functions $y$ are of class
$$C^{2r}([a,b], \mathbb{R}):=\left\{y:[a,
b] \cap \mathbb{T}\rightarrow\mathbb{R}\mid y^{\nabla^{2r}}\
\mbox{is continuous on}\ [a,b]_{\kappa^{2r}}\right\}.$$

\item $a,b \in \mathbb{T}$, $a<b$, and $[a,b]\cap \mathbb{T}$
has, at least, $2r+1$ points (\textrm{cf.} Remark~\ref{rem:nbr:points});

\item the Lagrangian $L(t,u_0,u_1,\ldots,u_r)$  has (standard)
    continuous partial derivatives with respect to
    $u_0,\ldots,u_r$, and partial nabla derivative with
    respect to $t$ of order $r+1$.

\end{enumerate}

\begin{remark}
For $\mathbb{T} = \mathbb{R}$ problem (\ref{problema
varacional}) coincides with the classical problem of the calculus of
variations with higher-order derivatives
(see, \textrm{e.g.}, \cite{Gelfand-Fomim:2000}).
Similar to $\mathbb{R}$, the results of the paper
are easily extended to the vectorial case, \textrm{i.e.},
to the case when admissible functions $y$ belong to
$C^{2r}([a,b], \mathbb{R}^n)$.
\end{remark}

\begin{remark}
For $r=1$, problem (\ref{problema varacional}) provides
the nabla analog of the delta problem on time scales
\eqref{eq:prb:cv:delta} studied in \cite{Zbig,Bohner:2004,zeidan:2004,Aveiro2Basia},
thus providing, in our opinion, a better formulation than that of
\cite{Atici:2006} (\textrm{cf.} \eqref{eq:prb:cv:nabla}).
For $r > 1$, restrictions to the time scale $\mathbb{T}$
must be done in order for (\ref{problema varacional})
to be well defined (\textrm{cf.} Remark~\ref{rem:rest:H} below).
\end{remark}

We begin with some technical results that will be useful in
the proof of our fundamental lemmas (\textrm{cf.} \S\ref{sub:sec:fl}).

\begin{proposition}
\label{teorema tecnico} Suppose that $a,b \in
\mathbb{T}$, $a < b$, and $f\in
C_{ld}([a,b], \mathbb{R}) $ is such that $f \geq 0$ on $[a,b]$. If
$\int_{a}^{b}f(t)\nabla t=0$, then $f=0$ on $[a,b]_\kappa$.
\end{proposition}

\begin{proof}
Suppose, by contradiction, that there exists $t_0\in[a,b]_\kappa$ such
that $f(t_0)>0$. If $t_0$ is left-scattered, then by the properties of the
integral (Theorems~\ref{propriedade integral} and
\ref{propriedades nabla integral}) we may conclude that
$$
\begin{array}{lcl}
\displaystyle{\int_{a}^{b}f(t) \nabla t} &= &
\displaystyle{\int_{a}^{\rho(t_0)}f(t)\nabla t +
\int_{\rho(t_0)}^{t_0}f(t)\nabla t + \int_{t_0}^{b}f(t)\nabla t}\\
& &\\ &\geq & \displaystyle{\int_{\rho(t_0)}^{t_0}f(t)\nabla
t=f(t_0)\left(t_0-\rho(t_0)\right)>0}
\end{array}
$$
which leads to a contradiction.
Suppose now that $t_0$ is left-dense. If $t_0 \neq a$, then by the
continuity of $f$ at $t_0$ we may conclude that there exists a
$\delta>0$ such that, for all $t\in]t_0- \delta,t_0]$, $f(t)>0$.
Since
$$\begin{array}{lcl}
\displaystyle{\int_{a}^{b}f(t) \nabla t} &=&
\displaystyle{\int_{a}^{t_0 -\delta}f(t) \nabla t} +
\displaystyle{\int_{t_0 -\delta}^{t_0}f(t) \nabla t} +
\displaystyle{\int_{t_0}^{b}f(t) \nabla t} \\ & &\\ &\geq &
\displaystyle{\int_{t_0-\delta}^{t_0}f(t)\nabla t >0}
\end{array}$$
we get a contradiction ($\delta$ may be chosen in such a way that $t_0
-\delta > a$). It remains to study the case when  $t_0 = a$. If
$t_0=a$ is right-dense, then by the continuity of $f$ at $t_0$ we
may conclude that there exists a $\delta>0$ such that, for all
$t\in[t_0,t_0 +  \delta [$, $f(t)>0$. Since
$$
\begin{array}{lcl}
\displaystyle{\int_{a}^{b}f(t) \nabla t} &=&
\displaystyle{\int_{t_0}^{t_0 +\delta}f(t) \nabla t} +
\displaystyle{\int_{t_0 +\delta}^{b}f(t) \nabla t} \\ & &\\ &\geq
& \displaystyle{\int_{t_0}^{t_0+\delta}f(t)\nabla t >0}
\end{array}
$$
we obtain again a contradiction.
Note that if $t_0=a$ and $a$ is
right-scattered, then $a\not\in[a,b]_{\kappa}$.
\end{proof}

\begin{remark}
In Proposition~\ref{teorema tecnico}
we cannot conclude that
$f=0$ on $[a,b]$. For example, consider $\mathbb{T}=\{1,2,3,4,5\}$
and $f(t)=1$ if $t=1$ and $f(t)=0$ otherwise. Clearly, $f$ is
continuous and $f\geq 0$ on $\mathbb{T}$. We have
$$
\int_{1}^{5} f(t)\nabla t = \sum_{t=2}^{5}f(t)=0 \, ,
$$
but $f\neq 0$ on $[1,5]$.
\end{remark}

From now on we restrict ourselves
to time scales $\mathbb{T}$ that satisfy the following condition $(H)$:
\begin{description}
\item[$(H)$] for each $t \in\mathbb{T}$, $(r-1) \left(\rho(t) - a_1t - a_0\right) = 0$ for some
$a_1\in\mathbb{R}^+$ and $a_0\in\mathbb{R}$.
\end{description}

\begin{remark}
\label{rem:rest:H}
Condition $(H)$ is equivalent to $r = 1$ or $\rho(t) = a_1t + a_0$
for some $a_1\in\mathbb{R}^+$ and $a_0\in\mathbb{R}$. Thus,
for the first order problem
of the calculus of variations we impose no restriction
on the time scale $\mathbb{T}$. For the higher-order problems
(\textrm{i.e.}, for $r > 1$) such restriction on the time scale is necessary.
Indeed, for $r > 1$ we are implicitly assuming in (\ref{problema varacional})
that $\rho$ is nabla differentiable, which is not true for a general
time scale $\mathbb{T}$.
 \end{remark}

\begin{remark}
\label{Remark 1}
Let $r > 1$.
Condition $(H)$ implies then that $\rho$ is nabla
differentiable. Hence,  $\nu$ is also nabla differentiable and
$\rho^{\nabla}(t)=a_1$, $t \in \mathbb{T}_\kappa$. Also note that
condition $(H)$ englobes the differential calculus
($\mathbb{T}=\mathbb{R}$, $a_1=1$, $a_0=0$), the difference
calculus ($\mathbb{T}=\mathbb{Z}$, $a_1=1$, $a_0=-1$), and the
q-calculus ($\mathbb{T}=\{ q^k: k \in
\mathbb{N}_0\}$ for some $q>1$, $a_1=\frac{1}{q}$, $a_0=0$).
\end{remark}

\begin{lemma}
\label{lemma5} Let $t\in \mathbb{T}_\kappa$ and $t\neq\max\mathbb{T}$
(if the maximum exists) satisfy the property
$\rho(t)<\sigma(t)=t$. Then, the backward jump operator $\rho$ is
not nabla differentiable at $t$.
\end{lemma}

\begin{proof}
We prove that $\rho$ is not continuous at $t\in \mathbb{T}_\kappa
\setminus\{\max\mathbb{T}\}$, which implies that $\rho$ is not
nabla differentiable at $t$. We begin by proving that
$\lim_{s\rightarrow t^+}\rho(s)=t$. Let $\varepsilon>0$ and take
$\delta=\varepsilon$. Then, for all $s\in ]t,t + \delta[$ we have
$|\rho(s)-t|\leq|s-t|<\delta=\varepsilon$. Since $\rho(t)\neq t$,
$\rho$ is not continuous at $t$.
\end{proof}

The following simply remark is very useful for our objectives:
\begin{remark}
\label{remark 0}
Since condition $(H)$ implies for $r > 1$ that $\rho$
is nabla differentiable, it follows from Lemma~\ref{lemma5}
that for the higher-order problem
$\mathbb{T}_\kappa \setminus\{\max\mathbb{T}\}$ cannot contain
points that are simultaneously right-dense and left-scattered.
\end{remark}

\begin{lemma}
\label{lemma6} Assume hypothesis $(H)$ and $r>1$. If
$f:\mathbb{T}\rightarrow \mathbb{R}$ is two times nabla
differentiable, then
$$
\label{Condition-lemma2} f^{\rho\nabla}(t)=a_1 f^{\nabla\rho}(t)\, ,
\quad t\in\mathbb{T}_{\kappa^2}.
$$
\end{lemma}

\begin{proof}
From Theorem~\ref{propriedades derivada} we know that
$f^{\rho}(t)=f(t)-\nu(t)f^\nabla(t)$. Thus,
$$f^{\rho\nabla}(t)=
\left(f(t)-\nu(t)f^\nabla(t)\right)^{\nabla}=
f^\nabla(t)-
\nu^{\nabla}(t)f^{\nabla\rho}(t)-\nu(t)f^{\nabla^2}(t)
=f^{\nabla\rho}(t)-\nu^{\nabla}(t)f^{\nabla\rho}(t).
$$
Since $\nu^{\nabla}(t)=1-a_1$, we may conclude that
$f^{\rho\nabla}(t)= a_1 f^{\nabla\rho}(t)$
for all $t\in\mathbb{T}_{\kappa^2}$.
\end{proof}

The next two lemmas are very useful for the proof of our higher-order
fundamental lemma of the calculus of variations on time scales
(Lemma~\ref{lemma7}).

\begin{lemma}\label{lemma funções admissíveis 1}
Assume that the time scale $\mathbb{T}$ satisfies condition $(H)$
and $\eta \in C^{2r}$ is such that
$\eta^{\nabla^{i}}(b)=0$ for all $i\in \{0,1,\ldots, r\}$.
Then, $\eta^{\rho \nabla^{i-1}}(b)=0$ for each $i\in \{1,\ldots, r\}$.
\end{lemma}

\begin{proof} If $b$ is left-dense, then the result is trivial
(just use Lemma~\ref{lemma6} and the fact that $\rho(b)=b$).
Suppose that $b$ is left-scattered and fix $i \in \{1,2,\ldots,
r\}$. From item 2 of Theorem~\ref{propriedades derivada}, we then
conclude that
$$\eta^{\nabla^{i}}(b)=
\left(\eta^{\nabla^{i-1}}\right)^{\nabla}(b)=
\frac{\eta^{\nabla^{i-1}}(b) -
\left(\eta^{\nabla^{i-1}}\right)^\rho(b)}{\nu(b)}.$$
Since
$\eta^{\nabla^{i}}(b)=0$ and $\eta^{\nabla^{i-1}}(b)=0$, then
$$
\left(\eta^{\nabla^{i-1}}\right)^\rho(b)=0.
$$
Lemma~\ref{lemma6} shows that
$$\left(\eta^{\nabla^{i-1}}\right)^\rho(b)=
\left(\frac{1}{a_1}\right)^{i-1}\left(\eta^\rho\right)^{\nabla^{i-1}}(b),$$
and we conclude that $\eta^{\rho\nabla^{i-1}}(b)=0$.
\end{proof}

\begin{lemma}\label{lemma funções admissíveis 2}
Assume that the time scale $\mathbb{T}$ satisfies condition $(H)$
and $\eta \in C^{2r}$ is such that
$$\eta^{\nabla^{i}}(\sigma^{r}(a))=0, \quad i\in \{0,1,\ldots,
r\}.$$
Then,
$$\eta^{\nabla^{i}}(\sigma^{i}(a))=0, \quad i\in \{0, 1,\ldots,
r-1\}.$$
\end{lemma}

\begin{proof}
If $a$ is right-dense, the result is trivial.
Suppose that $a$ is right-scattered (hence, $\sigma(a)$ is
left-scattered). Since $\rho$ is nabla differentiable,
by Remark~\ref{remark 0} we cannot have points that are simultaneously
right-dense and left-scattered. Hence, $\sigma(a)$, $\sigma^{2}(a)$,
$\sigma^{3}(a)$, $\ldots$, $\sigma^{r}(a)$ are left-scattered points.
By item \emph{2} of Theorem~\ref{propriedades derivada}, we
conclude that
$$\eta^{\nabla^{r}}(\sigma^{r}(a))=
\frac{\eta^{\nabla^{r-1}}(\sigma^{r}(a)) -
(\eta^{\nabla^{r-1}})^\rho(\sigma^{r}(a))}{\nu(\sigma^{r}(a))}=
\frac{\eta^{\nabla^{r-1}}(\sigma^{r}(a)) -
\eta^{\nabla^{r-1}}(\sigma^{r-1}(a))}{\nu(\sigma^{r}(a))}.$$
We have $\eta^{\nabla^{r}}(\sigma^{r}(a))=0$ and
$\eta^{\nabla^{r-1}}(\sigma^{r}(a))=0$. Then,
$$
\eta^{\nabla^{r-1}}(\sigma^{r-1}(a))=0.
$$
From item 2 of Theorem~\ref{propriedades derivada}
$$\eta^{\nabla^{r-1}}(\sigma^{r}(a))=
\frac{\eta^{\nabla^{r-2}}(\sigma^{r}(a)) -
(\eta^{\nabla^{r-2}})^\rho(\sigma^{r}(a))}{\nu(\sigma^{r}(a))}=
\frac{\eta^{\nabla^{r-2}}(\sigma^{r}(a)) -
\eta^{\nabla^{r-2}}(\sigma^{r-1}(a))}{\nu(\sigma^{r}(a))},$$
and using the hypothesis of the lemma we conclude that
$\eta^{\nabla^{r-2}}(\sigma^{r-1}(a))=0$. Since
$$\eta^{\nabla^{r-1}}(\sigma^{r-1}(a))=
\frac{\eta^{\nabla^{r-2}}(\sigma^{r-1}(a)) -
(\eta^{\nabla^{r-2}})^\rho(\sigma^{r-1}(a))}{\nu(\sigma^{r-1}(a))}=
\frac{\eta^{\nabla^{r-2}}(\sigma^{r-1}(a)) -
\eta^{\nabla^{r-2}}(\sigma^{r-2}(a))}{\nu(\sigma^{r-1}(a))},$$ we
obtain
$$
\eta^{\nabla^{r-2}}(\sigma^{r-2}(a))=0.
$$
Repeating recursively this process, we conclude the
proof.
\end{proof}


\subsection{Fundamental lemmas of the calculus of variations on time scales}
\label{sub:sec:fl}

We now present some fundamental lemmas of the
calculus of variations on time scales involving nabla derivatives.
This gives answer to a problem posed in \cite[\S3.2]{Tangier07}.
In what follows we assume that  $a$, $b \in \mathbb{T}$,
$a<b$, and $\mathbb{T}$ has sufficiently many points in order for
all the calculations to make sense.

\begin{lemma}
\label{lemma2}
Let $f \in C([a,b], \mathbb{R})$. If
$$
\label{formula2} \int_{a}^{b}f(t)\eta^{\nabla}(t)\nabla t=0 \ \ \
\ \ {\mbox for \ all} \ \ \eta \in C^1([a,b], \mathbb{R}) \ \
\mbox{ such \ that}\ \  \eta(a)=\eta(b)=0
$$
then
$$f(t)=c  \ \ \ \forall t\in [a,b]_\kappa$$
for some $c \in \mathbb{R}$.
\end{lemma}

\begin{proof}
Let $c$ be a constant defined by the condition
$$\int_{a}^{b}\left(f(\tau)-c\right)\nabla\tau=0$$
and let
$$\eta(t)=\int_{a}^{t}\left(f(\tau)-c\right)\nabla\tau.$$
Clearly, $\eta \in C^1([a,b], \mathbb{R})$ (by Theorem
\ref{antiderivada}, $\eta^{\nabla}(t)=f(t)-c$) and
$$\eta(a)=\int_{a}^{a}\left(f(\tau)-c\right)\nabla\tau=0 \ \ \ \
\mbox{and} \ \ \ \
\eta(b)=\int_{a}^{b}\left(f(\tau)-c\right)\nabla\tau=0.$$ Observe
that
$$\int_{a}^{b}\left(f(t)-c\right)\eta^{\nabla}(t)\nabla t=
\int_{a}^{b}\left(f(t)-c\right)^2\nabla t$$ and
$$\int_{a}^{b}\left(f(t)-c\right)\eta^{\nabla}(t)\nabla t=
\int_{a}^{b}f(t)\eta^{\nabla}(t)\nabla t- c
\int_{a}^{b}\eta^{\nabla}(t)\nabla t=0
-c\left(\eta(b)-\eta(a)\right)=0.$$ Hence,
$$\int_{a}^{b}\left(f(t)-c\right)^2\nabla t=0$$
which shows, by Proposition~\ref{teorema tecnico}, that
$$f(t)-c=0, \ \ \ \forall t \in [a,b]_\kappa.$$
\end{proof}

\begin{lemma}\label{lemma4}
Let $f, g \in C([a,b], \mathbb{R})$. If
$$
\label{formula3}  \int_{a}^{b}\left(f(t)\eta^\rho(t) + g(t)
\eta^{\nabla}(t)\right)\nabla t=0 \
$$
for all  $\eta \in C^1([a,b], \mathbb{R})$
such that $\eta(a)=\eta(b)=0$, then $g$  is nabla differentiable and
$$ g^{\nabla}(t)=f(t) \ \ \ \forall t\in [a,b]_\kappa.$$
\end{lemma}

\begin{proof} Define $A(t)=\int_{a}^{t} f(\tau)\nabla \tau$. Then
$A^{\nabla} (t)=f(t)$ for all $t \in [a,b]_\kappa$ (by Theorem
\ref{antiderivada}) and
$$\int_{a}^{b} A(t)\eta^{\nabla}(t) \nabla t=
\displaystyle \left[A(t)\eta(t)\right]_{t=a}^{t=b} -
\int_{a}^{b}A^{\nabla}(t) \eta^\rho (t) \nabla t= - \int_{a}^{b}
f(t) \eta^\rho(t) \nabla t$$ (by property \emph{8} of Theorem
\ref{propriedades nabla integral}). Hence,
$$
 \displaystyle{\int_{a}^{b}\left(f(t)\eta^\rho(t) + g(t)
\eta^{\nabla}(t)\right)\nabla t=0}  \Leftrightarrow \displaystyle{
\int_{a}^{b} \left(-A(t) + g(t) \right)\eta^{\nabla}(t)\nabla
t=0}.
$$
By Lemma \ref{lemma2} we may conclude that $-A(t) + g(t)=c$ for
all $t \in [a,b]_\kappa$ and  some $c \in \mathbb{R}$. Therefore,
$A^{\nabla}(t)=g^{\nabla}(t)$ for all $t \in [a,b]_\kappa$, proving the
desired result: $g^{\nabla}(t)=f(t)$ for all $ t\in [a,b]_\kappa.$
\end{proof}

\begin{remark}
If we consider $\mathbb{T}=\mathbb{R}$ in Lemmas~\ref{lemma2}
and \ref{lemma4}, we obtain some well known fundamental lemmas of
the classical calculus of variations (see, \textrm{e.g.},
\cite[pp.~10--11]{Gelfand-Fomim:2000}).
\end{remark}

\begin{remark}
Lemma~\ref{lemma2} remains true if $f$ is of class
$C_{ld}$ and the variation $\eta$ is of class
$C_{ld}^{1}$. Similar observation holds for
Lemma~\ref{lemma4}.
\end{remark}

We now prove a new and more  general fundamental lemma of the
calculus of variations. Lemma~\ref{lemma7} is used to prove our
Euler-Lagrange equation for variational problems on time scales
involving nabla derivatives of higher-order (Theorem
\ref{Euler-Lagrange for higher order}).

\begin{lemma}[higher-order fundamental lemma of the calculus of variations]
\label{lemma7}
Let $\mathbb{T}$ be a time scale satisfying condition $(H)$,
and $f_0$, $f_1$, $\ldots$, $f_r \in C([a,b], \mathbb{R})$. If
$$\int_{\sigma^{r-1}(a)}^{b}
\left(\sum_{i=0}^{r}f_i(t) \eta^{\rho^{r-i}\nabla^{i}}(t) \right)
\nabla t=0$$
for all $\eta \in C^{2r}([a, b], \mathbb{R})$ such that
\begin{align}
\eta\left(\sigma^{r-1}(a)\right)=0,& \quad \eta\left(b\right)=0,
\nonumber\\ &\vdots\nonumber \\
\eta^{\nabla^{r-1}}\left(\sigma^{r-1}(a)\right)=0,&\quad
\eta^{\nabla^{r-1}}\left(b\right)=0,\nonumber
\end{align}

then

$$\sum_{i=0}^{r} (-1)^i
\left(\frac{1}{a_1}\right)^{\frac{i(i-1)}{2}}f_i^{\nabla^i}(t) =0
\ ,  \ \ \ \ t \in [a,b]_{\kappa^{r}}.$$
\end{lemma}

\begin{proof}
We prove the lemma by mathematical induction.
If $r=1$, the result is true by Lemma~\ref{lemma4}.
Assume now that the result is true for some $r$, $r > 1$.
We want to prove that the result is then true for $r+1$.
Suppose that
$$\int_{\sigma^{r}(a)}^{b}
\left(\sum_{i=0}^{r+1}f_i(t) \eta^{\rho^{r+1-i}\nabla^{i}}(t)
\right) \nabla t=0$$ for all $\eta \in C^{2(r+1)}([a, b],
\mathbb{R})$ such that
$\eta\left(\sigma^{r}(a)\right)=0$, $\eta(b)=0$,
\ldots, $\eta^{\nabla^{r}}\left(\sigma^{r}(a)\right)=0$,
$\eta^{\nabla^{r}}(b)=0$.
We must prove that
$$\sum_{i=0}^{r+1} (-1)^i
\left(\frac{1}{a_1}\right)^{\frac{i(i-1)}{2}}f_i^{\nabla^i}(t) =0
\ ,  \ \ \ \
t \in [a,b]_{\kappa^{r+1}} \, .
$$
Note that
$$\int_{\sigma^{r}(a)}^{b}
\left(\sum_{i=0}^{r+1}f_i(t) \eta^{\rho^{r+1-i}\nabla^{i}}(t)
\right) \nabla t = \int_{\sigma^{r}(a)}^{b}
\left(\sum_{i=0}^{r}f_i(t) \eta^{\rho^{r+1-i}\nabla^{i}}(t)
 \right) \nabla t \ +  \int_{\sigma^{r}(a)}^{b}  f_{r+1}(t) \left
(\eta^{\nabla^r}
 \right)^{\nabla}(t)
\nabla t$$ and using the integration by parts formula (item 8 of
Theorem~\ref{propriedades nabla integral})
$$\int_{\sigma^{r}(a)}^{b}  f_{r+1}(t) \left (\eta^{\nabla^r}
\right)^{\nabla}(t)
\nabla t = \left[f_{r+1}(t)\eta^{\nabla^r}(t)
\right]^{t=b}_{t=\sigma^{r}(a)} - \int_{\sigma^{r}(a)}^{b}
f_{r+1}^{\nabla}(t)\left (\eta^{\nabla^r} \right)^{\rho}(t) \nabla
t.$$
Since $\eta^{\nabla^{r}}\left(\sigma^{r}(a)\right)=0$ and
$\eta^{\nabla^{r}}(b)=0$, we may conclude that
$$\int_{\sigma^{r}(a)}^{b}  f_{r+1}(t) \left (\eta^{\nabla^r}
\right)^{\nabla}(t)
\nabla t  = - \int_{\sigma^{r}(a)}^{b} f_{r+1}^{\nabla}(t)\left
(\eta^{\nabla^r} \right)^{\rho}(t) \nabla t. $$
By Lemma~\ref{lemma6},
$$\eta^{\nabla^{r} \rho}(t)= \left(
\frac{1}{a_1}\right)^{r}\eta^{\rho
\nabla^r}(t), \ \ t \in [a, b]_{\kappa^{r+1}}.$$ Hence,
$$\int_{\sigma^{r}(a)}^{b}  f_{r+1}(t) \left (\eta^{\nabla^r}
\right)^{\nabla}(t)
\nabla t  = - \int_{\sigma^{r}(a)}^{b} f_{r+1}^{\nabla}(t)\left(
\frac{1}{a_1}\right)^{r}\eta^{\rho \nabla^r}(t) \nabla t, $$
and
$$
\begin{array}{lcl}
 & \displaystyle \int_{\sigma^{r}(a)}^{b} & \left(
\displaystyle\sum_{i=0}^{r+1}f_i(t)
\eta^{\rho^{r+1-i}\nabla^{i}}(t) \right) \nabla t =\\
\\
& = & \displaystyle \int_{\sigma^{r}(a)}^{b}
\left(\sum_{i=0}^{r}f_i(t) \eta^{\rho^{r+1-i}\nabla^{i}}(t)
 \right) \nabla t -  \int_{\sigma^{r}(a)}^{b}
f_{r+1}^{\nabla}(t)\left(
\frac{1}{a_1}\right)^{r}\eta^{\rho \nabla^r}(t) \nabla t\\
\\
& = & \displaystyle\int_{\sigma^{r}(a)}^{b}
\left(\sum_{i=0}^{r-1}f_i(t) \left(\eta^\rho
\right)^{\rho^{r-i}\nabla^{i}}(t)
 + \left( f_r(t) -    f_{r+1}^{\nabla}(t)\left(
\frac{1}{a_1}\right)^{r}\right)   (\eta^\rho)^{\nabla^r}(t)
\right) \nabla t\\
\\
& = & \displaystyle \int_{\sigma^{r}(a)}^{\sigma^{r+1}(a)}
\left(\sum_{i=0}^{r-1}f_i(t) \left(\eta^\rho
\right)^{\rho^{r-i}\nabla^{i}}(t)
 + \left( f_r(t) -    f_{r+1}^{\nabla}(t)\left(
\frac{1}{a_1}\right)^{r} \right)     (\eta^\rho)^{\nabla^r}(t)
\right)    \nabla t\\
\\
 & \qquad + & \displaystyle\int_{\sigma^{r+1}(a)}^{b}
\left(\sum_{i=0}^{r-1}f_i(t) \left(\eta^\rho
\right)^{\rho^{r-i}\nabla^{i}}(t)
 + \left( f_r(t) -    f_{r+1}^{\nabla}(t)\left(
\frac{1}{a_1}\right)^{r}\right) (\eta^\rho)^{\nabla^r}(t)
 \right) \nabla t.
\end{array}
$$
We now prove that
\begin{equation}
\label{integral=0}\int_{\sigma^{r}(a)}^{\sigma^{r+1}(a)}
\left(\sum_{i=0}^{r-1}f_i(t) \left(\eta^\rho
\right)^{\rho^{r-i}\nabla^{i}}(t)
 + \left( f_r(t) -    f_{r+1}^{\nabla}(t)\left(
\frac{1}{a_1}\right)^{r}\right) (\eta^\rho)^{\nabla^r} (t) \right)
\nabla t \end{equation}
is equal to zero.
By Theorem~\ref{propriedade integral} the integral
(\ref{integral=0}) is equal to
\begin{multline*}
\left[\sum_{i=0}^{r-1}f_i(\sigma^{r+1}(a)) \left(\eta^\rho
\right)^{\rho^{r-i}\nabla^{i}}(\sigma^{r+1}(a))\right.\\
+ \left.\left( f_r(\sigma^{r+1}(a)) -
f_{r+1}^{\nabla}(\sigma^{r+1}(a))\left(
\frac{1}{a_1}\right)^{r}\right) (\eta^\rho)^{\nabla^r}
(\sigma^{r+1}(a))\right] \nu(\sigma^{r+1}(a)).
\end{multline*}
For each $i \in \{0, 1, \ldots, r\}$,
$$
\begin{array}{rcl}
\left(\eta^\rho \right)^{\rho^{r-i}\nabla^{i}}(\sigma^{r+1}(a)) &
=& \eta^{\rho^{r+1-i}\nabla^{i}}(\sigma^{r+1}(a))\\ & = &
(a_1)^{i(r+1-i)} \eta^{\nabla^{i} \rho^{r+1-i}}(\sigma^{r+1}(a)) \
\ \ \ \ \mbox{(by Lemma~\ref{lemma6})}\\ &=& (a_1)^{i(r+1-i)}
\eta^{\nabla^{i} }(\sigma^{i}(a))\\ &=& 0 \ \ \ \ \
\mbox{(by Lemma~\ref{lemma funções admissíveis 2})}
\end{array}
$$
proving that the integral (\ref{integral=0}) is equal to zero.
Then,
$$
\begin{array}{lcl}
 & \displaystyle \int_{\sigma^{r}(a)}^{b} & \displaystyle
\left(\sum_{i=0}^{r+1}f_i(t)
\eta^{\rho^{r+1-i}\nabla^{i}}(t) \right) \nabla t \\
\\
&= & \displaystyle \int_{\sigma^{r+1}(a)}^{b}
\left(\sum_{i=0}^{r-1}f_i(t) \left(\eta^\rho
\right)^{\rho^{r-i}\nabla^{i}}(t)
 + \left( f_r(t) -    f_{r+1}^{\nabla}(t)\left(
\frac{1}{a_1}\right)^{r}\right) (\eta^\rho)^{\nabla^r}(t) \right)
\nabla t.
\end{array}
$$
Observe that,
$$
\begin{array}{rcl}
\eta^\rho\left(\sigma^{r+1}(a)\right)& = & \eta(\sigma^r(a))=  0
\\ & & \\ (\eta^\rho)^{\nabla}\left(\sigma^{r+1}(a)\right)& = &
a_1\eta^{\nabla}(\sigma^r(a))=  0 \\ &\vdots& \nonumber \\
(\eta^\rho)^{\nabla^{r-1}}\left(\sigma^{r+1}(a)\right) & = &
(a_1)^{r-1}\eta^{\nabla^{r-1}}(\sigma^r(a))=0\\
\end{array}
$$
and, by Lemma \ref{lemma funções admissíveis 1},
$$
\begin{array}{rcl}
\eta^\rho(b)& =& 0 \\ & & \\ (\eta^\rho)^{\nabla}(b)& = & 0 \\
&\vdots& \nonumber \\ (\eta^\rho)^{\nabla^{r-1}}(b)& = & 0.\\
\end{array}
$$
Then, by the induction hypothesis, we conclude that
$$\sum_{i=0}^{r-1} (-1)^i
\left(\frac{1}{a_1}\right)^{\frac{i(i-1)}{2}}f_i^{\nabla^i}(t) \ +
\ (-1)^r \left(\frac{1}{a_1}\right)^{\frac{r(r-1)}{2}}
\left(f_r(t) -
f^{\nabla}_{r+1}(t)\left(\frac{1}{a_1}\right)^r\right)^{\nabla^r}
(t) =0 \ , \ \ \ \  t \in [a,b]_{\kappa^{r+1}} \, ,$$
which is equivalent to
$$\sum_{i=0}^{r+1} (-1)^i
\left(\frac{1}{a_1}\right)^{\frac{i(i-1)}{2}}f_i^{\nabla^i}(t)
=0 \ , \ \ \ \  t \in [a,b]_{\kappa^{r+1}}.$$
\end{proof}

\begin{remark}
The differentiability of
functions $f_0$, $f_1$, \ldots, $f_r$ is not assumed a priori.
\end{remark}


\subsection{Euler-Lagrange equations for higher-order problems}
\label{sub:sec:EL:ho}

Before presenting the Euler-Lagrange equation for the variational
problem \eqref{problema varacional},
we introduce  the following definition.

\begin{definition} We say that $y_{\ast}\in C^{2r}([a,b],
\mathbb{R})$ is a weak local minimizer (respectively weak local
maximizer) for problem \eqref{problema varacional} if there exists
$\delta
>0$ such that
$$
\mathcal{L}[y_{\ast}]\leq \mathcal{L}[y] \ \ \ \ \ \ \
(\mbox{respectively} \ \ \  \mathcal{L}[y_{\ast}]\geq
\mathcal{L}[y])
$$
for all $y \in C^{2r}([a,b], \mathbb{R})$ satisfying the boundary
conditions in \eqref{problema varacional} and
$$
\parallel y - y_{\ast}\parallel_{r,\infty} < \delta \, ,
$$
where
$$
\parallel y\parallel_{r,\infty}:= \sum_{i=0}^{r}
\parallel y^{\rho^{r-i}\nabla^{i}}\parallel_{\infty}
$$
and
$$ \parallel y\parallel_{\infty} :=\sup_{t \in
[a,b]_{\kappa^r}}\mid y(t) \mid.
$$
\end{definition}

\begin{remark}
\label{rem:nbr:points}
Observe that if $[a,b]$ has $2r$ points, \textrm{i.e.},
$$[a,b]=\left\{\rho^{2r-1}(b), \rho^{2r}(b), \ldots,
\rho^{2}(b),\rho(b), b\right\},$$ then
$$
\begin{array}{rcl}
\mathcal{L}[y(\cdot)] &= & \displaystyle
\int_{\sigma^{r-1}(a)}^{b}
L(t,y^{\rho^r}(t),y^{\rho^{r-1}\nabla}(t),\ldots,
y^{\rho\nabla^{r-1}}(t),y^{\nabla^r}(t))\nabla t \\ & & \\ &= &
\displaystyle \int_{\rho^{r}(b)}^{b}
L(t,y^{\rho^r}(t),y^{\rho^{r-1}\nabla}(t),\ldots,
y^{\rho\nabla^{r-1}}(t),y^{\nabla^r}(t))\nabla t \\ & & \\ &=&
\displaystyle \sum_{i=0}^{r-1} \displaystyle
\int_{\rho^{i+1}(b)}^{\rho^{i}(b)}L(t,y^{\rho^r}(t),y^{\rho^{r-1}\nabla}(t),\ldots,
y^{\rho\nabla^{r-1}}(t),y^{\nabla^r}(t))\nabla t\\ & & \\ &=&
\displaystyle \sum_{i=0}^{r-1} L(\rho^{i}(b),
y^{\rho^r}(\rho^{i}(b)), y^{\rho^{r-1}\nabla}(\rho^{i}(b)),
\ldots, y^{\nabla^r}(\rho^{i}(b)))\left( \rho^{i}(b) -
\rho^{i+1}(b) \right)\\ & & \\ &=& \displaystyle \sum_{i=0}^{r-1}
L(\rho^{i}(b), y(\rho^{r+i}(b)),
(a_1)^{r-1}y^{\nabla}(\rho^{r+i-1}(b)), \ldots,
y^{\nabla^r}(\rho^{i}(b)))\left( \rho^{i}(b) - \rho^{i+1}(b)
\right).
\end{array}
$$
Using the boundary conditions
in \eqref{problema varacional} and the formula
$$y^{\nabla}(t)=\frac{y(t)-y(\rho(t))}{t-\rho(t)} \, ,$$
it is then possible to calculate
$$y(\rho^{r+i}(b)), y^{\nabla}(\rho^{r+i-1}(b)), \ldots,
y^{\nabla^r}(\rho^{i}(b))$$
for all $i=0,1, \ldots, r-1$. Therefore,  the above integral is
constant for every admissible function $y(\cdot)$. We conclude that if
$[a,b]$ has only $2r$ points, the problem is trivial (because there is
nothing to minimize or maximize). For this reason, we are assuming
that $[a,b]$ has, at least, $2r+1$ points.
\end{remark}

We are in conditions to prove the following
first-order necessary optimality condition for problems of the
calculus of variations with higher-order nabla derivatives:

\begin{theorem}[Euler-Lagrange equation for problem (\ref{problema varacional})]
\label{Euler-Lagrange for higher order} Let $\mathbb{T}$ be a time
scale satisfying hypothesis $(H)$, and $a,b \in \mathbb{T}$,
$a<b$, with $[a,b] \cap \mathbb{T}$ containing, at least, $2r+1$
points. If $y_\ast$ is a weak local extremum (minimizer or
maximizer) of problem (\ref{problema varacional}), then $y_\ast$
satisfies the Euler-Lagrange equation
$$
\sum_{i=0}^{r}(-1)^i\left(\frac{1}{a_1}\right)^{\frac{i(i
-1)}{2}}L^{\nabla^i}_{u_i}\left(t,y_\ast^{\rho^r}(t),y_\ast^{\rho^{r-1}\nabla}(t),
\ldots,y_\ast^{\rho\nabla^{r-1}}(t),y_\ast^{\nabla^r}(t)\right)=0 \, ,
$$
$t\in[a,b]_{\kappa^{2r}}$.
\end{theorem}

\begin{proof}Suppose that $y_\ast$ is a weak
local minimizer (resp. maximizer) for problem (\ref{problema varacional}).
Let $\eta\in C^{2r}([a,b], \mathbb{R})$ be an admissible variation, \textrm{i.e.},
$\eta$ is such that $\eta, \eta^{\nabla}, \ldots ,
\eta^{\nabla^{r-1}}$ vanish at $t=\sigma^{r-1}(a)$ and $t=b$.
Defining $\phi:\mathbb{R} \rightarrow \mathbb{R}$ by
$\phi(\epsilon):= \mathcal{L}[y_\ast + \epsilon \eta]$ it is clear
that $\phi$ has a minimum (resp. maximum) at $\epsilon=0$ and,
therefore,
$$\phi^\prime(0)=0.$$
Since
$$\phi(\epsilon)=\int_{\sigma^{r-1}(a)}^{b}
L(t,y_{\ast}^{\rho^r}(t)+ \epsilon
\eta^{\rho^r}(t),y_{\ast}^{\rho^{r-1}\nabla}(t)+ \epsilon
\eta^{\rho^{r-1}\nabla}(t),\ldots,y_{\ast}^{\nabla^r}(t) +
\epsilon \eta^{\nabla^r}(t))\nabla t \, ,$$
then
$$\phi^\prime(0)=0 \Leftrightarrow
\int_{\sigma^{r-1}(a)}^{b} \left(\sum_{i=0}^{r}L_{u_i}(\bullet)
\eta^{\rho^{r-i}\nabla^{i}}(t) \right) \nabla t=0$$
where $L_{u_i}$ denote the partial derivative of $L(t, u_0, u_1,
\ldots, u_r)$ with respect to $u_i$ and we write, for brevity,
$$(\bullet)=
(t,y_{\ast}^{\rho^r}(t),y_{\ast}^{\rho^{r-1}\nabla}(t),\ldots,
y_{\ast}^{\rho\nabla^{r-1}}(t),y_{\ast}^{\nabla^r}(t)).$$
By Lemma~\ref{lemma7}, we conclude that
$$\sum_{i=0}^{r} (-1)^i
\left(\frac{1}{a_1}\right)^{\frac{i(i-1)}{2}}L_{u_i}^{\nabla^i}(\bullet)
=0 \ ,
\ \ \ \ t \in \left([a,b]_{\kappa^{r}}\right)_{\kappa^{r}}=[a,b]_{\kappa^{2r}} \, ,
$$
proving the intended result.
\end{proof}

As a straight corollary to Theorem~\ref{Euler-Lagrange for higher order},
we give the Euler-Lagrange equation for the higher-order
variational problem of q-calculus:

\begin{corollary}[the q-calculus Euler-Lagrange equation]
\label{sub:sec:example} Fix $q>1$, $\mathbb{T}= \{ q^k: k \in
\mathbb{N}_0\}$. Let $a, b \in \mathbb{T}$ such that $[a,b]$
contains at least $2r+1$ points. If $y_{\ast}$ is a weak local
extremum for the correspondent problem \eqref{problema
varacional}, then $y_{\ast}$ satisfies the Euler-Lagrange equation
\begin{equation*}
\sum_{i=0}^{r}(-1)^i q^{\frac{i(i
-1)}{2}}L^{\nabla^i}_{u_i}\left(t,y_\ast(q^{-r}t),
q^{1-r}y_\ast^{\nabla}(q^{1-r} t), \ldots,
q^{-1}y_\ast^{\nabla^{r-1}}(q^{-1}t),y_\ast^{\nabla^r}(t)\right)=0
\end{equation*}
for all $t\in [q^{2r}a,b]$.
\end{corollary}



\end{document}